\documentclass[11pt]{article}

\usepackage[backref=page]{hyperref}  %

\usepackage{amsmath}  
\usepackage{amsthm} 
\usepackage{amssymb} 
\usepackage{hyperref} 
\usepackage{color} 

\newtheorem{thm}{Theorem}
\newtheorem{inspr}[thm]{}

\newenvironment{definitie}{\begin{itemize}\item[ ]\hspace{-26pt}\bf Definition \rm }{\end{itemize}}
\newenvironment{notatie}{\begin{itemize}\item[ ]\hspace{-26pt}\bf Notation \rm }{\end{itemize}}
\newenvironment{voorbeeld}{\begin{itemize}\item[ ]\hspace{-26pt}\bf Example \rm }{\end{itemize}}
\newenvironment{stelling}{\begin{itemize}\item[ ]\hspace{-26pt}\bf Theorem \rm }{\end{itemize}}
\newenvironment{propositie}{\begin{itemize}\item[ ]\hspace{-26pt}\bf Proposition \rm }{\end{itemize}}
\newenvironment{lemma}{\begin{itemize}\item[ ]\hspace{-26pt}\bf Lemma \rm }{\end{itemize}}
\newenvironment{opmerking}{\begin{itemize}\item[ ]\hspace{-26pt}\bf Remark \rm }{\end{itemize}}
\newenvironment{voorwaarde}{\begin{itemize}\item[ ]\hspace{-26pt}\bf Condition \rm }{\end{itemize}}

\newcommand{\defin}{\begin{inspr}\begin{definitie}}  
\newcommand{\edefin}{\end{definitie}\end{inspr}}

\newcommand{\notat}{\begin{inspr}\begin{notatie}}  
\newcommand{\enotat}{\end{notatie}\end{inspr}}

\newcommand{\voorb}{\begin{inspr}\begin{voorbeeld}}  
\newcommand{\evoorb}{\end{voorbeeld}\end{inspr}}

\newcommand{\stel}{\begin{inspr}\begin{stelling}}
\newcommand{\estel}{\end{stelling}\end{inspr}}

\newcommand{\prop}{\begin{inspr}\begin{propositie}}
\newcommand{\eprop}{\end{propositie}\end{inspr}}

\newcommand{\lem}{\begin{inspr}\begin{lemma}}
\newcommand{\elem}{\end{lemma}\end{inspr}}

\newcommand{\opm}{\begin{inspr}\begin{opmerking}}
\newcommand{\eopm}{\end{opmerking}\end{inspr}}

\newcommand{\voorw}{\begin{inspr}\begin{voorwaarde}}
\newcommand{\evoorw}{\end{voorwaarde}\end{inspr}}

\newcommand{\bew}{\vspace{-0.3cm}\begin{itemize}\item[ ] \bf Proof\rm: }
\newcommand{\ebew}{\hfill $\qed$ \end{itemize}}

\newcommand{\ssnl}{\vskip 3pt} 
\newcommand{\snl}{\vskip 7pt} 
\newcommand{\nl}{\vskip 12pt} 

\newcommand{\ot}{\otimes}

\newcommand{\tussenen}{\qquad\quad\text{and}\qquad\quad}

\newcommand{\rood}{\color{red}}
\newcommand{\blauw}{\color{blue}}
\newcommand{\zwart}{\color{black}}

\numberwithin{thm}{section}   
\numberwithin{equation}{section} 

\usepackage{stmaryrd}
\usepackage[textsize=scriptsize]{todonotes}
\usepackage[all]{xy}
\definecolor{teal}{rgb}{0.0, 0.5, 0.5}

\def\Hom{{\sf Hom}}
\def\End{{\sf End}}
\def\ol{\overline}
\def\joost{\color{teal}}

\newcommand{\mycomment}[2]{{\blauw {#1}\ssnl}{\zwart{#2}}} 
\newcommand{\keepcomment}[1]{}
\newcommand{\oldcomment}[1]{}

\begin{document}

\setcounter{section}{-1}  

%
 %

\centerline{\bf \Large Multiplier algebras and local units}
\vspace{13pt}
\centerline{\it Alfons Van Daele \rm $^{(*)}$  and \it Joost Vercruysse \rm $^{(**)}$}
\bigskip\bigskip
{\bf Abstract} 
\nl
Let $A$ be an algebra over any field.
We do not assume that $A$ has an identity. The \emph{multiplier algebra} $M(A)$ is a unital algebra associated to $A$.
If we require the product in $A$ to be non-degenerate (as a bilinear form), the multiplier algebra can be characterized as the largest algebra containing $A$ as an essential ideal. We recall the basic definitions and provide some more information about this notion. 
\ssnl
We endow the multiplier algebra $M(A)$ with the {\it strict topology}. Then we show that $A$ is dense in $M(A)$ if and only if there exist local units in $A$.
\ssnl
We include various examples. In particular, we are interested in the underlying algebras of multiplier Hopf algebras, algebraic quantum groups, algebraic quantum hypergroups, weak multiplier Hopf algebras and algebraic quantum groupoids. In all these cases, one can show that the algebras have local units. We have also included some examples arising from co-Frobenius coalgebras.
\ssnl
For most of the material treated in this note, it is only the ring structure of the algebra that plays a role. For this reason, we develop the theory here for rings. But they are not required to have an identity for the multiplicative structure.
\ssnl

\nl 
Date: {\it July 11, 2025}
\nl

\vskip 6cm
\nl
\hrule
\medskip
\begin{itemize}
\item[$^{(*)}$] Department of Mathematics, KU Leuven, Celestijnenlaan 200B,
B-3001 Heverlee, Belgium. E-mail: \texttt{alfons.vandaele@kuleuven.be}
\item[$^{(**)}$] Département de Mathématiques, Université Libre de Bruxelles, Boulevard du Triomphe, B-1050 Bruxelles, Belgium.
E-mail: \texttt{joost.vercruysse@ulb.be}
\end{itemize}

\newpage

%
%

\section{\hspace{-17pt}. Introduction} \label{s:introduction}  

Rings without unit element appear naturally in all areas of Mathematics. Often add hoc arguments are put in place to treat the lack of a unit element in the considered setting. This has made that similar results, sometimes with different terminology have been used in various fields, and also that some relevant results well-known in some areas remained unknown to others. The aim of the present paper is twofold. Firstly, we want to provide a comprehensive overview of some basic properties and techniques of two central notions for non-unital rings: their multiplier algebra and the existence of local units. Secondly, we provide some new results and connections, filling certain gaps in the existing literature on the subject. We illustrate the theory with examples arising from coalgebra and multiplier Hopf algebra theory.

Let $A$ be an associative algebra over a field $k$, not necessarily having a unit element but such that the product in $A$ is non-degenerate. We can always define the multiplier algebra $M(A)$.  Because the algebra is assumed to be non-degenerate, it can be characterized as the biggest unital algebra that contains $A$ as an essential ideal. We recall the notion in Definition \ref{defin:nd} and some of its properties in Section \ref{s:multipliers}.
\ssnl
The multiplier algebra $M(A)$ can be endowed with the strict topology. 
For this topology,  a net $(x_\alpha)$ of elements in $M(A)$ converges to an element $x\in M(A)$  if for all $a\in A$ we have $x_\alpha a=xa$ and $ax_\alpha =ax$ when $\alpha$ is large enough. It is not hard to see that $A$ is strictly dense in $M(A)$ if and only if $A$ has local units. This result is well-known in the context of operator algebras. See e.g.\ Section 3.12 in \cite{P}, in particular  item 3.12.17 (Notes and remarks), with a reference to \cite{Bu}.
\ssnl
Recall that $A$ is said to have local units if for any finite number of elements $a_1,a_2,\dots,a_n$ in $A$ there exists an element $e\in A$ so that $a_i e=e a_i =a$ for all $i$. Observe that the existence of local units implies that the multiplication is non-degenerate. If $A$ has local units, then it also must be an idempotent algebra which means that any element of $A$ is the sum of products of elements of $A$. The condition is written as $A=A^2$. Moreover, $A$ is also firm in the sense that the multiplication induces an isomorphism $A\ot_A A\cong A$.
\ssnl
If there are no local units, we still can consider the closure $\overline A$ of $A$ in $M(A)$ for the strict topology. It will be a proper ideal of $M(A)$, still containing  $A$ as an ideal. All of this is treated in Section \ref{s:local-units} of the paper. 
\ssnl
Basic and motivating examples are included in these two sections. However, in a separate section, Section \ref{s:examples}, we give examples arising from semiperfect coalgebras, from the theory of multiplier Hopf algebras and algebraic quantum groups (i.e.\ multiplier Hopf algebras with integrals), algebraic quantum hypergroups and of weak multiplier Hopf algebras and algebraic quantum groupoids. In fact, this note finds its origin in the study of (weak) multiplier Hopf algebras. 
\ssnl
In another, related paper \cite{VD-Ve}, we study infinite matrix algebras. We describe again the multiplier algebra for these examples. But in that paper, we discuss another aspect of non-unital algebras. For finite-dimensional algebras, we have the notion of a Frobenius algebra. This means that there exists a faithful functional. The situation for non-unital, possibly infinite-dimensional algebras is more complex. This is treated in \cite{VD-Ve} and illustrated with these infinite matrix algebras. Infinite matrix algebras appear naturally in the theory of algebraic quantum groups and as examples in the theory of weak multiplier Hopf algebras.
\ssnl
Most of the theory we discuss in this paper works for rings and not only for algebras. Therefore, we will introduce the concepts and prove the result in this more general setting. It makes no difference when we are dealing with algebras and one can safely work with the underlying ring structure. 
\ssnl
This paper, as well as the one on infinite matrix algebras \cite{VD-Ve} does not contain many new results. The purpose is rather to collect some of the known results, scattered in the literature, and explain some of the features in a more systematic way and by examples.
\nl
\bf Content of the paper \rm
\nl
As we mentioned already, in Section \ref{s:multipliers} we begin with reviewing the definitions of the  multiplier algebras $M(A)$, $L(A)$ and $R(A)$ for a non-degenerate algebra. We recall the extension properties of a non-degenerate homomorphism from an algebra $A$ to the multiplier algebra $M(B)$ of another algebra $B$. We also deal with some special cases. We give a few basic examples throughout the exposition but at the end of the section, we include a more specific (finite-dimensional) example that allows us to illustrate some special features. 
\ssnl
In fact, as only the ring structure is important, we do this for rings and not only for algebras.
\ssnl
In Section \ref{s:local-units} we start with the precise definition of local units and we recall some sufficient condition for the existence of local units.  The main part is devoted to the introduction of the strict topology on the multiplier algebra $M(A)$ of an algebra $A$. We show that $A$ is dense in $M(A)$ if and only if $A$ has local units. Also here we illustrate the result with some basic examples. Also in this section, we do all this for rings and not only for algebras.
\ssnl
In Section \ref{s:examples} we focus on the underlying algebras of multiplier Hopf algebras, algebraic quantum hypergroups and weak multiplier Hopf algebras. It has been proven at various places in the literature that local units exist for these algebras. We collect some of the given arguments, scattered in the literature and discuss the ideas behind many of these proofs. We start this section however with some examples of a different nature, namely the rational dual of a co-Frobenius coalgebra. 
\snl
Finally, in Section \ref{s:conclusions}, we draw  conclusions and discuss some open problems for future research. We also refer once more to related work. 
\nl
\bf Notations and conventions, basic references \rm
\nl
In this paper, we  discuss the material for general rings and not only for (associative) algebras over any field $k$. Recall however that in many earlier papers on multiplier Hopf algebras and similar objects, only algebras over $\mathbb C$  were considered.
\ssnl
Sometimes, we also discuss involutive rings. In that case, we assume that the ring carries an involution $a\mapsto a^*$. It is an additive map satisfying $(a^*)^*=a$ and $(ab)^*=b^*a^*$ for all $a,b\in A$. 
\ssnl
The rings and algebras need not be unital. But the multiplication is {\it always assumed to be non-degenerate} (as a bilinear form). This is automatic if there is a unit. It is also automatic if there are  local units. 
\ssnl
Often our rings will be idempotent, in the sense that any element is a sum of products of elements. Again this is automatic if it is unital, or more general, when it has local units. This condition will not be imposed, but it is often a consequence of the other conditions that are considered.
\ssnl
When $P$ and $Q$ are subspaces of a ring, we denote by $PQ$ the subspace spanned by products $pq$ with $p\in P$ and $q\in G$. With this convention, a ring $A$ is idempotent if $A^2=A$.
\ssnl
We use $M(A)$ for the multiplier ring of $A$. We recall the notion in the beginning of Section \ref{s:introduction}. The identity in $M(A)$ is always denoted by $1$ while we use $\iota$ for the identity map. 
\ssnl
The multiplier algebra of a non-degenerate algebra, as we use it here, is considered in \cite{VD-mha} but it should be mentioned that it has been studied earlier (see e.g.\ \cite{Da}). For the notion of a comultiplication, as it appears in the theory of coalgebras, we refer to \cite{A}, \cite{S} and \cite{R}. For the notion of a  comultiplication in the theory of multiplier Hopf algebra, we refer to \cite{VD-mha} and in the setting for weak multiplier Hopf algebras, to \cite{VD-W0} and \cite{VD-W1}. \oldcomment{Ik heb nog niet systematisch het woord coproduct vervangen door comultiplication zoals je gesuggereerd heb.}{}
\ssnl
By a {\em homomorphism} between two rings $A$ and $B$, we mean an additive and multiplicative map. When $A$ and $B$ are unital rings, then we call a homomorphism {\em unital} if it sends the unit of $A$ to the unit of $B$.
\nl
\bf Acknowlegdements \rm
\nl
The first named author, being retired for many years, likes to thank the KU Leuven for the opportunity to continue doing research at the mathematical institute. Further, the two authors enjoy the cooperation among each other. 

%
%

\section{\hspace{-17pt}. Multiplier rings and  algebras} \label{s:multipliers} 

Let $A$ be a ring. We \emph{do not assume} that it has an identity (for the multiplication).  On the other hand we do not want the product to be completely trivial either. In principle, it could happen for a ring $A$ that $ab=0$ for all $a$ and all $b$ in $A$ and such a situation is not very interesting. Therefore, in this paper we will assume that all our rings and algebras have a {\em non-degenerate multiplication} (except if we explicitly mention the contrary). We recall the definition.

\defin\label{defin:nd}
We call a ring $A$ {\it non-degenerate} if its multiplication is non-degenerate as a bilinear form. This means that, for a given element $a\in A$, we have $a=0$ if either $ab=0$ for all $b\in A$ or $ba=0$ for all $b\in A$. We call $A$ {\em left non-degenerate} if $ab=0$ for all $b\in A$ implies $a=0$ and {\em right non-degenerate} if $ba=0$ for all $b\in A$ implies $a=0$.
\edefin

If the ring $A$ has an identity, the condition is trivially satisfied. So all rings with identity are non-degenerate in the above sense. The same is true if the ring has local units. We will treat this case in Section \ref{s:local-units}, were we also recall the precise definitions. 
\ssnl
Here are a few simple examples. 

\voorb\label{voorb:1.2}
i) Let $X$ be any set. Consider the algebra $K(X)$ of complex functions with finite support in $X$, with pointwise operations. This is a non-degenerate algebra with local units. It has a unit if and only if $X$ is finite.
\vskip 3pt
ii) Take for $X$ the set $\mathbb N$ of natural numbers. Then elements in $K(X)$ are sequences $(a_n)_{n=1}^\infty$ with $a_n\in \mathbb C$ but so that only finitely many of these elements are non-zero. We can also consider 
the bigger algebra $c_0$ of all sequences that tend to $0$ when $n\to\infty$. This is also a non-degenerate algebra, but now there are no local units. 
\evoorb

A ring is non-degnerate exactly if it is both left and right non-degenerate. It is not too difficult to find an algebra that is left non-degenerate, but not right non-degenerate,
see e.g.\ Remark 2.5 in  \cite{VD-Ve}. Hence, for a ring to be non-degenerate, we really need the condition on the two sides. \oldcomment{\rood We moeten op het einde alle interne en externe referenties nakijken.}{}
\ssnl
Many properties, well-known in the case of finite-dimensional unital algebras are no longer true for possibly non-unital or infinite-dimensional algebras. Also some of the notions have to be reformulated. 
\nl
\bf The multiplier ring $M(A)$ of $A$ \rm
\nl
The notion of a multiplier in the case of a non-unital algebra with a non-degenerate multiplication is introduced in the first paper on multiplier Hopf algebras (\cite{VD-mha}). The treatment was inspired by the concept as used in the theory of operator algebras (see e.g. Section 3.12 in \cite{P}). 
In fact, the notion had been introduced already before in 1969  by J.\ Dauns in \cite{Da}.
\ssnl
Here we construct the multiplier ring $M(A)$ of a non-degenerate ring $A$  as given (for algebras)  in the appendix of \cite{VD-mha}).

\notat\label{notat:1.3}
Let $A$ be a non-degenerate ring. Denote by $M(A)$ the set of  pairs $(\rho,\lambda)$ of maps from $A$ to itself satisfying $\rho(a)b=a\lambda(b)$ for all $a,b\in A$. 
We can define a map $A\to M(A)$ by $a\mapsto (\rho_a,\lambda_a)$ where
\begin{equation*}
\rho_a(b)=ba
\tussenen
\lambda_a(b)=ab
\end{equation*}
for all $b$. Clearly, because the multiplication is assumed to be non-degenerate, this map is injective. We will identify $A$ with its image in $M(A)$ under this embedding.
\ssnl
When $x=(\rho,\lambda)$ is any element in $M(A)$, we call $x$ a {\em multiplier} for $A$ and we write 
$$xa := \lambda(a) \tussenen ax := \rho(a).$$ 
The defining relation of a multiplier can then be written more intuitively as an associativity condition
$(ax)b=a(xb)$ when $x\in M(A)$ and $a,b\in A$. 
\enotat
The following result is important (and easy to prove).

\lem\label{lem:1.4a}
Assume that $x$ is a multiplier of a non-degenerate ring $A$. For any $a,a',b,b'\in A$ and $x\in M(A)$, we have
\begin{equation*}
(a+a')x = ax+a'x\tussenen x(b+b')=xb+xb' 
\end{equation*}
and
\begin{equation*}
(aa')x=a(a'x) \tussenen x(bb')=(xb)b'.
\end{equation*}
In other words, for a multiplier $x=(\lambda,\rho)$, we see that $\lambda$ is a right $A$-linear endomorphsm and $\rho$ is a left $A$-linear endomorphism. 
\ssnl
If $A$ is an algebra over a field $k$, the maps $a\mapsto ax$ and $a\mapsto xa$ are also $k$-linear.
\elem
\bew
i) For any $a,a',b$ in $A$ we have
\begin{equation*}
a(x(b+b'))=(ax)(b+b')=(ax)b+(ax)b'=a(xb)+a(xb')=a(xb+xb')
\end{equation*}
and from the non-degeneracy of the multiplication, we get $x(b+b')=xb+xb'$. In a similar way one proves that $(a+a')x = ax+a'x$.
\ssnl
ii) 
For any $a,a',b\in A$ and $x\in M(A)$, using the above conventions, we find
$$((aa')x)b= (aa')(xb)= a(a'(xb))=a((a'x)b)=(a(a'x))b$$
and because the multiplication is non-degenerate, we find
$(aa')x=a(a'x)$.
Similarly, one obtains that
$x(bb')=(xb)b'$
whenever $x\in M(A)$ and $a,b,b'\in A$.
\ebew

Remark that 
$M(A)$ can be defined, also if the ring does not have a non-degenerate multiplication, but then we can not embed $A$ in $M(A)$.  That case is not very useful.
\ssnl
As the above notation suggests, we would like to think about $\lambda(a)=xa$ as a multiplication of the element $x\in M(A)$ with the element $a$ of $A$. The following observation shows that this interpretation is justified, because this multiplication is associative.
\ssnl
Also remark that a multiplier $x$ is defined by saying what $xa$ and $ax$ are for $a\in A$ and by then verifying that $(ax)b=a(xb)$ for all $a,b\in A$.
\ssnl
With this point of view, it is not hard to verify the following results. 
\prop
Let $A$ be a non-degenerate ring and consider $M(A)$ as above. Then we have the following properties:
\ssnl
i) $M(A)$ is made into a unital associative algebra if we define
\begin{equation*}
(xy)a=x(ya)
\qquad\qquad\text{and}\qquad\qquad
a(xy)=(ax)y
\end{equation*}
whenever $x,y\in M(A)$ and $a\in A$. The identity in $M(A)$ is given exactly by the identity maps on $A$.
\ssnl
ii) $A$ is a two-sided ideal in $M(A)$, which is moreover an {\em essential} ideal in the sense that for $x\in M(A)$, we have $x=0$ if either $xa=0$ for all $a$ in $A$ or if $ax=0$ for all $a$ in $A$.
\ssnl
iii) In the case that $A$ is a ring with an involution $a\mapsto a^*$ , then also $M(A)$ carries an involutive structure, given by
\begin{equation*}
x^*a=(a^*x)^*
\tussenen
ax^*=(xa^*)^*
\end{equation*}
for $x\in M(A)$ and $a\in A$.
\eprop

\defin
We call $M(A)$ the \emph{multiplier ring} of $A$. If $A$ is an algebra, we call it the \emph{multiplier algebra}.
\edefin

Remark that Lemma \ref{lem:1.4a} implies in particular that the ring underlying the multiplier algebra of an algebra $A$ coincides with the multiplier ring of the underlying ring of $A$. Henceforth, we can use the notation $M(A)$ unambiguously, and use the terms multiplier algebra and multiplier ring interchangeably. 
\ssnl
We now show that $M(A)$ is characterized by the following universal property.

\prop\label{prop:Misbiggest}
For any unital ring $B$ that contains $A$ as an ideal, there exists a unique unital homomorphism $j:B\to M(A)$ whose restriction to $A$ is the natural embedding of $A$ in $M(A)$.  In other words, the homomorphism $j$ makes the following diagram commutative.
\[
\xymatrix{
A \ \ar@{^(->}[dr] \ar@{^(->}[rr] && M(A) \\
& B \ar@{.>}[ur]_j
}
\]
Moreover, $j$ is injective if and only if $A$ is essential as an ideal in $B$ and consequently, $M(A)$ is the largest unital algebra that contains $A$ as an essential two-sided ideal.
\eprop

\bew
Suppose that $B$ is a unital algebra, containing $A$ as a two-sided ideal. Define a map $j:B\to M(A)$ by
\begin{equation*}
 a(j(b))=ab
\qquad\qquad\text{and}\qquad\qquad
(j(b))a=ba
\end{equation*}
whenever $a\in A$. It is clear that $j(b)\in M(A)$ as 
$$(a(j(b)))a'=(ab)a'=a(ba')=a((j(b))a')$$
 for all $a,a'$. One can verify that $j$ is multiplicative and that $j(1_B)=1_{M(A)}$. 
By definition, the injectivity of $j$ means that $A$ is an essential ideal. The last statement is then an immediate consequence.
 \ebew

We can also look at the ring of left multipliers and the ring of right multipliers. Within the spirit of the notations above, we introduce them as follows. 

\defin 
A {\em left multiplier} $y$ for $A$ is an additive map $a\mapsto ya$ from $A$ to itself satisfying $y(ab)=(ya)b$ for all $a,b\in A$. In other words, a left multiplier for $A$ is right $A$-linear endomorphism of $A$. We denote by $L(A)$ the set $\End_A(A)$ of all left multipliers for $A$. Similarly, a {\em right multiplier} $z$ for $A$ is a left $A$-linear endomorphism of $A$, that is, $z$ is an addive map $a\mapsto az$ from $A$ to itself satisfying $(ab)z=a(bz)$ for all $a,b\in A$. The set ${_A\End}(A)$ of all right multipliers for $A$ is denoted as $R(A)$.
\edefin

As we have seen in Lemma \ref{lem:1.4a}, for a multiplier $x$, the linear map $a\mapsto xa$ is a left multiplier in the sense of the above definition, while $a\mapsto ax$ is a right multiplier. 
\ssnl
Also remark that, just as in the case of multipliers, we can define left and right multipliers also when the product is degenerate. We will not consider this case further.
\ssnl
The usual composition of maps induces a (unital) ring structure on both $L(A)$ and $R(A)$. The obvious morphisms 
$$A\to L(A),\ a\mapsto (b\mapsto ab), \qquad A\to R(A),\ a\mapsto (b\mapsto ba)$$
are multiplicative and moreover injective if $A$ left, respectively right, non-degenerate.
\ssnl

Intuitively, a multiplier for $A$ is at the same time a left and right multiplier. Moreover a given multiplier $x=(\lambda,\rho)$ is completely determined by its left multiplier $\lambda$, or equally by its right multiplier $\rho$. Hence one can think about $M(A)$ as the ``intersection'' of $L(A)$ and $R(A)$. However, this intersection is not just the intersection of $L(A)$ and $R(A)$ viewed as subsets of the set of all additive endomorphisms of $A$. Indeed, this intersection corresponds to the set of all $A$-bilinear endomorphisms of $A$, which, in case of a unital algebra $A$, is naturally isomorphic to the center of $A$. In order to make the characterization of $M(A)$ as an intersection of $L(A)$ and $R(A)$ precise, let us recall the following result from \cite[Section 1.2]{JaVe}. We will illustrate this phenomenon in an example at the end of this section.

\prop
Let $A$ be a non-degenerate ring. Then we have injective maps
\begin{eqnarray*}
\overline{(-)}:&&R(A)\to {_A\Hom_A}(A\ot A,A),\ \overline x(a\ot b)=(ax)b\\
\underline{(-)}:&&L(A)\to {_A\Hom_A}(A\ot A,A),\ \overline y(a\ot b)=a(yb)
\end{eqnarray*}
The space $M(A)$ of multipliers of $A$ is then given by the pullback (as abelian group) of the above morphisms.
\[\xymatrix{
M(A) \ar[rr] \ar[d] && R(A) \ar[d]^-{\overline{(-)}} \\
L(A) \ar[rr]_-{\underline{(-)}} && {_A\Hom_A}(A\ot A,A)
}
\]
\eprop
\bew
Let us prove that $\overline{(-)}$ is injective. Suppose $x\in R(A)$ is such that $\overline x=0$. Then we find for all $a,b\in A$ that $(ax)b=0$. By non-degeneracy of the multiplication we then find that $ax=0$ for all $a$, hence $x=0$ in $R(A)$. In the same way, we find that $\underline{(-)}$ is injective. The rest of the statement follows directly from the definitions. 
\ebew

If we consider again the  simple Example \ref{voorb:1.2}, we easily find the following result.

\voorb\label{voorb:1.2bis}
Let $X$ be an infinite set and $A$ the algebra $K(X)$ of complex functions with finite support in $X$. The multiplier algebra can be identified with the algebra $C(X)$ of all complex functions on $X$.
\evoorb

If we take for $X$ a locally compact space and for $A$ the algebra of continuous complex functions with compact support, then the multiplier algebra is the algebra of all continous complex functions on $X$.
\ssnl
The following result is however less obvious.

\voorb\label{voorb:2.11c}
Let $A$ be the algebra $c_0$ of sequences, tending to $0$ as we have in item ii) of Example \ref{voorb:1.2}. It can be shown that the multiplier algebra here is $\ell^\infty$, the algebra of bounded sequences. It is clear that the pointwise product $ab\in c_0$ when $a$ is a bounded sequence and $b\in c_0$. To prove the converse, first observe that  any multiplier of $c_0$ must be given by a sequence. Further it is not so difficult to give an example where $a$ is not bounded, $b$ is in $c_0$ but $ab$ is not in $c_0$. This proves the converse, namely that any multiplier of $c_0$ is given by an element of $\ell^\infty$.
\evoorb

Similarly, the multiplier algebra of $C_0(X)$ of continuous complex functions on a locally compact space $X$, tending to $0$ at infinity, is the algebra of bounded continuous complex function on $X$.

We consider this example again in Section \ref{s:local-units}, see Example \ref{voorb:2.3}. More examples can be found in forthcoming \cite{VD-Ve}

\nl
\bf Non-degenerate homomorphisms \rm
\nl
Recall the following notion, formulated for algebras in the appendix in \cite{VD-mha}. 

\defin \label{defin:1.10}
Let $A$ and $B$ be non-degenerate rings. Assume that $\gamma:A\to M(B)$ is a homomorphism. It is called non-degenerate if $\gamma(A)B=B$ and $B\gamma(A)=B$
\edefin

If $A$ and $B$ are unital rings, then $\gamma:A\to M(B)=B$ is non-degenerate if and only if it sends the unit of $A$ to the unit of $B$. This shows that non-degenerate homomorphisms as in Definition \ref{defin:1.10} provide a correct analogue of unital homomorphisms in the non-unital setting.
\ssnl
It was observed in \cite{JaVe} that non-degenerate homomorphisms $\gamma:A\to M(B)$ can also be defined without making reference to the multiplier algebra. Let us recall this now.

\lem\label{lem:JV}
Let $A$ and $B$ be non-degenerate rings. Then there is a bijective correspondence between non-degenerate homomorphisms $\gamma:A\to M(B)$ and $A$-bimodule structures on $B$, such that the multiplication of $B$ is $A$-bilinear and $A$-balanced, and the bimodule structure is {\em idempotent} in the sense that $AB=B=BA$.
\elem
\bew
Given $\gamma$, we define an $A$-bimodule structure on $B$ by 
$$a.b=\gamma(a)b\tussenen b.a=b\gamma(a).$$
 Then the definition and properties of multipliers tell us that the multiplication of $A$ is indeed $A$-bilinear and $A$-balanced. The non-degeneracy of $\gamma$ corresponds exactly with the idempotency of the bimodule structure. 

Conversely, given an $A$-bimodule structure on $B$, we can define $\gamma(a)b=a.b$ and $b\gamma(a)=b.a$. 
\ebew

With this notion at hand, we have the following result (see Proposition A.5 in \cite{VD-mha} and Theorem 1.15 in \cite{JaVe}).
\ssnl 
Let us first introduce some notation. By ${\sf NdI}$ we denote the category whose objects are non-degenerate idempotent rings (recall that a ring $A$ is idempotent if $A^2=A$) and where morphisms are the non-degenerate homomorphisms as in Definition \ref{defin:1.10}. Let ${\sf Rng}$ be the category of unital rings with unital homomorphisms between them.

\prop\label{prop:1.9}
i) If $\gamma:A\to M(B)$ is a non-degenerate homomorphism, then there is  a unique unital 
homomorphism $\gamma_1: M(A) \to M(B)$ such that $\gamma$ coincides with the composition of $\gamma_1$ with the embedding $A\to M(A)$. In other words, $\gamma_1$ makes the following diagram commute
\[
\xymatrix{
A \ar[rr]^-\gamma \ar@{^(->}[dr] && M(B)\\
& M(A)\ar@{.>}[ur]_-{\gamma_1}
}
\]
ii) There is a functor $M:{\sf NdI}\to {\sf Rng}$, which sends a non-degneratre idempotent ring to its multiplier ring and a homomorphism to its extension as in part i) of this proposition. Moreover, $M$ is a right adjoint to the obvious forgetful functor ${\sf Rng}\to {\sf NdI}$.
\eprop

\bew
i) Consider a multiplier $m\in M(A)$. By the non-degeneracy of $\gamma$, it follows that if $\gamma_1$ exists, then $\gamma_1(m)\in M(B)$ should be defined by means of the formulas
\begin{eqnarray*}
\gamma_1(m)b&=& \sum \gamma_1(m)\gamma(a_i)b_i=\sum \gamma(ma_i)b_i \\
d\gamma_1(m)&=&\sum d_i\gamma(c_i)\gamma_1(m)=\sum d_i\gamma(c_im)
\end{eqnarray*}
where $b = \sum \gamma(a_i)b_i$ and $d= \sum d_i\gamma(c_i)$. We use that $B=\gamma(A)B=B\gamma(A)$. 
\ssnl
To see that $\gamma_1$ defined as above is well-defined,
first assume that we have elements $a_i\in A$ and $b_i$ in $B$ such that $\sum \gamma(a_i)b_i=0$. Then for all 
$c\in A$ and $d\in B$ we have
\begin{equation*}
d\gamma(c)\sum_i \gamma(ma_i)b_i=\sum_i d\gamma(cma_i)b_i=d\gamma(cm)\sum_i \gamma(a_i)b_i=0.
\end{equation*}
Because $B\gamma(A)=B$ and because the product in $B$ is non-degenerate, it follows that also $\sum \gamma(ma_i)b_i=0$. 
\ssnl
Applying a symmetric argument for the right multiplier $\gamma_1(m)$, we find that the formulas above are independent from the chosen representatives in $\gamma(A)B$ or $B\gamma(A)$. 
\ssnl
We clearly have 
\begin{equation*}
d\gamma(c)\gamma(ma)b=d\gamma(cm)\gamma(a)b
\end{equation*}
and so, $\gamma_1(m)$ is indeed a multiplier.
\ssnl
It is now straightforward to complete the proof of the first statement.
\ssnl
ii) If $A$ is an idempotent ring, then the embedding from $A$ in $M(A)$ is a non-degenerate homomorphism in the sense of Definition \ref{defin:1.10}. This embedding is the unit morphism in the category ${\sf NdI}$. Its extension as in item i) of the proposition is then clearly the identity map from $M(A)$ to itself. One can also verify that the construction from item i) preserves the composition and henceforth we obtain the stated functor. Furthermore, for a unital ring $R$ and a non-degenerate ring $A$, we find that non-degenerate homomorphisms $R\to M(A)$ are exactly unital homomorphisms, which is exactly the desired adjunction property.
\ebew

It is a common practice to denote this extension $\gamma_1$ again by the symbol $\gamma$. We will do this in what follows.
\ssnl
In case $A$ is not idempotent, that is $A^2\neq A$,
we can still extend the embedding $A\to M(A)$ to the identity map from $M(A)$ to itself, but this may not be the only extension. 
\ssnl
Remark that in \cite{JaVe}, the adjunction property between $M$ and the forgetful functor was erroneously stated the other way around.

\voorb
Consider two sets $X$ and $Y$ and a map $\alpha:X\to Y$. Let $A$ be the algebra $K(Y)$ and $B$ the algebra $K(X)$. Define $\gamma(f)(x)=f(\alpha(x))$ for $x\in X$ and $f\in A$. It could happen that the function $\gamma(f)$ no longer has finite support. But surely, we get a homomorphism from $A$ to $M(B)$.
\ssnl
One can observe that $\delta_x=\delta_x\gamma(\delta_{\alpha(x)})$ for all $x\in X$. Moreover, any $f\in B$ having finite support can be written as a finite sum $f=\sum f\delta_x$, where $x$ varies over the support of $f$. Then we see that $f=\sum f\delta_x\gamma(\delta_{\alpha(x)})\in B\gamma(A)$ and henceforth the homomorphism from $A$ to $M(B)$ is non-degenerate.
\ssnl
The extension here is obviously given by the same formula $\gamma(f)(x)=f(\alpha(x))$ for $x\in X$, but now $f$ is allowed to be any function on $Y$.
\evoorb

We are interested in the following special case.

\prop\label{prop:1.11}
Assume that $A$ is a subalgebra of a unital algebra $B$ such that $AB=B$ and $BA=B$. Then $A$ is automatically non-degenerate and its multiplier algebra can be realized as a subalgebra of $B$. In fact we have that an element $b$ of $B$ belongs to $M(A)$ if and only if $ba\in A$ and $ab\in A$ for all $a\in A$.
\eprop

\bew
i) First we show that $A$ is non-degenerate. For this, let $a\in A$ and assume $ac=0$ for all $c\in A$. Then $acd=0$ for all $c\in A$ and all $d\in B$. As $AB=B$ we get $ab=0$ for all $b\in B$. Then $a=0$. Similarly on the other side.
\vskip 3pt
ii) Use $\gamma$ to denote the embedding of $A$ in $B$. Here $M(B)=B$ because $B$ is unital. The conditions $AB=B$ and $BA=B$ precisely mean that $\gamma$ is non-degenerate. By the previous result, we can extend $\gamma$ to the multiplier algebra $M(A)$ of $A$.
\vskip 3pt
iii) We claim that this extension is still injective. Indeed, assume $m\in M(A)$ and $\gamma(m)=0$. Then $\gamma(ma)b=0$ for all $a\in A$ and all $b\in B$. So $\gamma(ma)=0$ and $ma=0$ for all $a\in A$. It follows that $m=0$.
\vskip 3pt
iv) Finally, let $b\in B$. If $b\in M(A)$, then $ba\in A$ and $ab\in A$. Conversely, if $ba\in A$ and $ab\in A$, we can define a multiplier $m\in M(A)$ by $ma=ba$ and $ab=am$. Obviously, we will get $m=b$.
\ebew

Let us remark that the non-degeneracy condition $AB=B=BA$ is rather strong. In case $A$ is an algebra with local units, it implies that $A$ has a unit, and is a unital subalgebra of $B$. Indeed, let us write $1=\sum a_ib_i\in AB$ and let $e\in A$ be a left local unit for the $a_i$ appearing in this sum. Then we find that $e=e.1=e.\sum a_ib_i=\sum ea_ib_i=\sum a_ib_i=1$, so $1\in A$. In particular, the proposition cannot be applied to an algebra $A$ of the form $K(X)$, where $X$ is an infinite set. The latter can also be seen directly as follows. Let $K(X)\subset B$ for a unital algebra $B$ and $1=\sum a_ib_i$, where $a_i\in K(X)$ and $b_i\in B$. Each of the $a_i$ has finite support. Denote $Y\subset X$ the union of all these supports, then $Y$ is still a finite set. If $X$ is infinite, there exists $x\in X\setminus Y$. Then $\delta_x a_i=0$ for all $i$ and hence $0=\sum \delta_xa_ib_i=\delta_x 1 =\delta_x$, which is a contradiction. 
\ssnl
However, Proposition \ref{prop:1.11} applies to general non-degenerate algebras and we will provide a non-trivial illustrative example at the end of this section.
\ssnl
We have similar results for one-sided multiplier algebras.

\prop\label{prop:1.19a}
 Let $A$ and $B$ be non-degenerate algebras. And let $\gamma:A\to L(B)$ be a homomorphism of $A$ in the algebra of left multipliers on $B$. Then $B$ can be considered as a left $A$-module via the action $a.b=\gamma(a)b$ (see also Lemma \ref{lem:JV}). 
Assume now that $B$ is firm as left $A$-module, that is, the action of $A$ on $B$ induces an isomorphism $A\ot_A B\cong B$ (this is for example the case if $A$ has local units and $\gamma(A)B=B$). Then we can extend $\gamma$ uniquely to a {\em unital} homomorphism from $L(A)$ to $L(B)$.
\eprop

\bew
For any $x\in L(A)$, we define $\gamma(x)\in L(B)$ as the following composition
$$\xymatrix{
B \cong A\ot_A B \ar[rr]^-{x\ot_A \iota_B} && A\ot_A B \cong B
}$$
Remark that since $x$ is a left multiplier, it is a right $A$-linear map, so that $x\ot_A \iota_B$ is indeed well-defined. Explicitly, we have that 
$$\gamma(x)b=\sum (xa_i)b_i$$
whenever $b$ is the image of $\sum a_i\ot_A b_i$ under the isomorpism $A\ot_A B\cong B$.
\ebew

By similar arguments as in the proof for the two-sided case Proposition~\ref{prop:1.11}, the previous proposition can be used to characterize the left multiplier algebra of $A$ as a subalgebra of $B$, if $B$ has a unit.

\prop\label{prop:1.19b}
Assume that $A$ is a subalgebra of a unital algebra $B$ such that $B$ is firm as left $A$-module, $A\ot_A B\cong B$. Then the left multiplier algebra $L(A)$ embeds in $B$. And an element $b$ of $B$ is in (the image of) $L(A)$ if an only if $ba\in A$ for all $a\in A$.
\eprop

The firmness property is not very easy to check in practice. Therefore, we provide a sufficient condition for this to hold which is more practical to verify.

\prop\label{prop:1.19}
Consider a unital algebra $B$ and a (possibly non-unital) subalgebra $A\subset B$. Denote $B^L=\{b\in B ~|~ ba\in A, \forall a\in A\}$, which is a unital subalgebra of $B$, and suppose that $AB^L=B^L$. Then $B$ is firm as left $A$-module and therefore $B^L\cong L(A)$.
\eprop

\bew
One can easily observe that $B^L$ is indeed a unital subalgebra of $B$. Now write $1=\sum a_ib_i$, with $a_i\in A$ and $b_i\in B^L$. We claim that the map $B\to A\ot_A B, b\mapsto \sum_i a_i\ot_A b_ib$, is an inverse of the morphism $A\ot_A B\to B,\ a\ot_A b\mapsto ab$, showing that $B$ is firm as left $A$-module. Indeed for any $b\in B$, we have $\sum_i a_ib_ib= 1b=b$. Conversely, for any $a'\ot_A b'\in A\ot_A B$, since $b_ia'\in A$, we find $\sum_i a_i \ot_A b_i a'b' = \sum_i a_ib_ia'\ot_A b'=a'\ot_A b'$.
\ssnl
We can therefore apply Proposition~\ref{prop:1.19b}, to conclude that $B^L\cong L(A)$.
\ebew

\snl
\bf An example \rm
\nl
In our paper {\it Infinite Matrix Algebras} \cite{VD-Ve}, we consider a certain type of non-degenerate algebras and describe in detail the multiplier algebras, together with some other features.
\ssnl
In this paper, we just include the following special example.

\voorb\label{voorb:1.16}
Let $M_3$ be the algebra  of $3$ by  $3$ complex matrices. Consider a set of matrix units
\begin{equation*}
\{e_{ij}\mid i,j=1,2,3\}.
\end{equation*}
Recall that $e_{ij}e_{k\ell}=\delta_{jk}e_{i\ell}$ where $\delta_{jk}=0$ if $j\neq k$ and $\delta_{jj}=1$.
\ssnl
i) Now let $A$ be the subspace of $M_3$ spanned by the matrix units
\begin{equation*}
e_{11},e_{12},e_{13},e_{23},e_{33}.
\end{equation*}
It is the space of upper triangular matrices, however with $0$ in the middle of the diagonal. It is not hard to see that this is a subalgebra of $M_3$. Indeed
\begin{equation*}
\left(\begin{matrix} a & b & c \\ 0 & 0 & d \\ 0 & 0 & e \end{matrix}\right)
\left(\begin{matrix} a' & b' & c' \\ 0 & 0 & d' \\ 0 & 0 & e' \end{matrix}\right)
=
\left(\begin{matrix} aa' & ab' & ac'+bd'+ce' \\ 0 & 0 & de' \\ 0 & 0 & ee' \end{matrix}\right).
\end{equation*}
ii) This algebra is {\it idempotent} as we have
\begin{align*}
e_{11}&=e_{11}e_{11}
\tussenen
e_{12}=e_{11}e_{12}\\
e_{13}&=e_{11}e_{13}=e_{12}e_{23}=e_{13}e_{33}\\
e_{23}&=e_{23}e_{33}
\tussenen
e_{33}=e_{33}e_{33}.
\end{align*}
iii) The algebra is {\it non-degenerate}. To show this, take $x\in A$ and write
\begin{equation*}
x=ae_{11}+be_{12}+ce_{13}+de_{23}+ee_{33}.
\end{equation*}
Now assume that $yx=0$ for all $y\in A$. Then 
\begin{align*}
e_{11}x&=a e_{11} +be_{12}+ce_{13}=0\\
e_{12}x&= de_{13}=0\\
e_{33}x&= ee_{33}=0. 
\end{align*}
It follows that all co\"efficients are $0$ and so $x=0$.  On the other hand, assume that $xy=0$ for all $y\in A$. Then
\begin{align*}
xe_{11}&=a e_{11}=0\\
xe_{23}&= be_{13}=0\\
xe_{33}&= ce_{13}+de_{23}+ee_{33}=0. 
\end{align*}
And again follows that $x=0$.
\ssnl
iv) The algebra has {\it no local units}. Let again $x=ae_{11}+be_{12}+ce_{13}+de_{23}+ee_{33}$. Then $xe_{23}=be_{13}$ and there is no element $x$ in the algebra satisfying $xe_{23}=e_{23}$. 
\ssnl
In particular, the algebra has no unit. Since it is finite dimensional, it would have local units if and only if it would have a unit.
\ssnl
v) We claim that the {\it left multipliers} are spanned by the algebra itself and the matrices $e_{22}$ and $e_{32}$. 
\ssnl
First, it is easy to verify that these are left multipliers. On the other hand, assume that $\lambda$ is a linear map from $A$ to itself satisfying $\lambda(xy)=\lambda(x)y$ for all $x,y\in A$. Then 
\begin{equation*}
\lambda(e_{11})=\lambda(e_{11}e_{11})=\lambda(e_{11})e_{11}
\end{equation*}
and it follows that there is a scalar $p$ so that $\lambda(e_{11})=pe_{11}$. Then 
\begin{equation*}
\lambda(e_{12})=\lambda(e_{11}e_{12})=\lambda(e_{11})e_{12}=pe_{11}e_{12}
\end{equation*}
and so also 
$\lambda(e_{12})=pe_{12}$.
Similarly 
\begin{equation*}
\lambda(e_{13})=\lambda(e_{11}e_{13})=\lambda(e_{11})e_{13}=pe_{11}e_{13}
\end{equation*}
and $\lambda(e_{13})=pe_{13}$. Now
\begin{align*}
\lambda(e_{13})&=\lambda(e_{12}e_{23})=\lambda(e_{12})e_{23}=pe_{12}e_{23}=pe_{13}\\
\lambda(e_{13})&=\lambda(e_{13}e_{33})=\lambda(e_{13})e_{33}=pe_{13}e_{33}=pe_{13}
\end{align*}
but this gives no furhter information. It just confirms earlier results. On the other hand
\begin{equation*}
\lambda(e_{23})=\lambda(e_{23}e_{33})=\lambda(e_{23})e_{33}
\end{equation*}
and this implies that $\lambda(e_{23})=ue_{13}+ve_{23}+ we_{33}$ for some $u,v,w\in \mathbb C$.
Finally
\begin{equation*}
\lambda(e_{33})=\lambda(e_{33}e_{33})=\lambda(e_{33})e_{33}
\end{equation*}
and also this implies that $\lambda(e_{23})=u'e_{13}+v'e_{23}+ w'e_{33}$ for some $u',v',w'\in \mathbb C$.
The above shows that the space of left multipliers is at most $7$-dimensional. Since we already know the space of left multipliers contains $A$ ($5$-dimensional) and also the extra matrix units $e_{22}$ and $e_{32}$, we know that these generate all.
\ssnl
vi) A completely similar argument will give that the algebra of right multipliers is spanned by $A$ and the matrices $e_{21}$ and $e_{22}$
\ssnl
vi) We have an example where there is a left multiplier that is not a right multiplier and vice versa. The {\it multiplier algebra} here is spanned by $A$ and the matrix $e_{22}$. 
\ssnl
vii) Writing $B=M_3$, one can check that the conditions of Proposition \ref{prop:1.11} and Proposition \ref{prop:1.19} are satisfied. First observe that $B^L$ is the algebra generated by the elements of $A$ and matrix units $e_{22}$ and $e_{32}$. As before, one easily checks that the product of any of these elements with an element of $A$ is again in $A$, hence is contained in $B^L$. On the other hand, we find that $e_{21}e_{12}=e_{22}\not\in A$ and $e_{31}e_{12}=e_{32}\not\in A$, so no other elements of $B$ belong to $B^L$.
Similarly, one observes that $\{b\in B~|~ ab, ba \in A, \forall a\in A\}$ is spanned by $A$ and $e_{22}$. 
Furthermore, we have that $AB=B$ and $AB^L=B^L$ since 
\begin{eqnarray*}
e_{21}&=&e_{23}e_{31} \in A B\\
e_{22}&=&e_{23}e_{32} \in A B^L\\
e_{31}&=&e_{33}e_{31} \in A B\\
e_{32}&=&e_{33}e_{32} \in A B^L
\end{eqnarray*}
Similarly, one checks that $BA=B$. From this we can deduce by Proposition \ref{prop:1.11} that $M(A)$ can be identified with the algebra of upper triangular matrices, and by Proposition \ref{prop:1.19} that $B^L$ can be identified with $L(A)$, confirming the observations made above.
\evoorb

Let us describe the extension- and restriction-of-scalars functors, associated to the homomorphism $A\to M(A)$, which provides a pair of adjoint functors between the category of (possibly non-unital) modules over $A$ and the category of unital modules over its multiplier algebra.

\prop
Let $A$ be a (non-degenerate) algebra and $M(A)$ its multiplier algebra. Denote by ${_A\sf Mod}$ the category of (possibly non-unital) left $A$-modules and by ${_{M(A)}\sf Mod}^u$ the category of unital left $M(A)$-modules. \\
i) Then the extension-of-scalars functor
$$M(A)\ot_A -: {_A \sf Mod}\to {_{M(A)} \sf Mod}^u$$
is a left adjoint to the restriction-of-scalars functor ${_{M(A)} \sf Mod}^u\to {_A \sf Mod}$. 
\ssnl 
ii) If $A$ is idempotent, then $A\cong M(A)\ot_A A$ via the unit of the adjunction.
\ssnl 
iii) If $A$ is idempotent, then for any firm left $A$-module $M$, the unit of the adjunction is an isomorphism $M\cong M(A)\ot_A M$.
It follows that firm modules over an idempotent algebra form a full subcategory of the category of unital $M(A)$-modules. 
\eprop

\bew
i) The adjunction can be deduced from the following natural isomorphism for any left $A$-module $M$ and unital left $M(A)$-module $N$:
\begin{eqnarray*}
{_{M(A)} \sf Hom}(M(A)\ot_A M,N) & \cong & {_A \sf Hom}(M,N)\\
\varphi & \mapsto & \varphi(1\ot_A -)\\
(x\ot m \mapsto xf(m)) & \mapsfrom & f
\end{eqnarray*}
The unit and counit of adjunction are then given by 
\begin{eqnarray*}
\eta_M: M\to M(A)\ot_A M,&& \eta_M(m)= 1\ot_A m\\
\epsilon_N: M(A)\ot_A N\to N,&& \epsilon_N(x\ot n)=xn
\end{eqnarray*}
\ssnl
ii) Define $\ol\eta_A: M(A)\ot_A A\to A, \ol\eta_A(x\ot a)=xa$. We claim that this is a two-sided inverse for $\eta_A$. For any $a\in A$, we find $\ol\eta_A\circ\eta(a)=1a=a$. On the other hand, using idempotency of $A$, write $a=\sum a_ia'_i$. Then for any $x\ot_A a\in M(A)\ot_A A$, we obtain
\begin{eqnarray*}
\eta_A\circ\ol\eta_A(x\ot_A a) &=& 1\ot_A xa = \sum 1\ot_A xa_ia'_i \\
&=& \sum xa_i\ot_A a'_i = x\ot \sum a_ia'_i=x\ot a
\end{eqnarray*}
and hence $\eta_A$ is an isomorphism.
\ssnl
iii) Using part ii), together with the firmness of $M$, we find an isomorphism
$$M\cong A\ot_A M \cong M(A)\ot_A A\ot_A M\cong M(A)\ot_A M,$$
which coincides with the unit $\eta_M$.
The last statement then expresses that $M(A)\ot_A-$ is a fully faithful left adjoint.
\ebew

%
%

\section{\hspace{-17pt}. Local units and multiplier algebras}\label{s:local-units} 

In this section, we consider a non-degenerate ring $A$ together with its multiplier ring $M(A)$ as we have reviewed in the previous section.
\ssnl
We will now focus  on the existence of local units and relate the property with the {\it density} of $A$ in $M(A)$. This is now in a topological sense, explained in an item  below. As mentioned already in the introduction, the result is inspired by a corresponding result for operator algebras (cfr.\ Section 3.12 in \cite{P}).
\nl
\bf Local units for $A$ \rm
\nl

There exist various definitions of what is called a \emph{ring with local units}. In each of the cases, one requires for any finite number of elements in the ring, the existence of an (other) element acting as the identity on the chosen elements. However, sometimes one requires these `local units' to posses stronger properties, such as being idempotents \cite{AM}, central or commuting idempotents \cite{Ab}, or the existence of a complete set of orthogonal idempotents. The definition we state here, is equivalent to the notion of what is called an s-unital ring in, for example, \cite{T}. For an overview of various notions of rings with local units and more references to literature, we refer to \cite{Ve} and \cite{N}

\defin\label{defin:2.1}
We say that the ring $A$ has {\it local units} if for any finite subset $a_1, a_2, \dots, a_n$ of $A$ there is an element $e\in A$ so that 
\begin{equation*}
 a_ie=a_i \qquad\qquad\text{and}\qquad\qquad ea_i=a_i
\end{equation*}
for all $i=1,2,\dots,n$.
\edefin

It is well-known (see e.g. \cite[Theorem 1]{T}, or \cite[Lemma 2.2]{Ve} for a more general version)
that a weaker condition is sufficient to guarantee the existence of local units as above. We not only formulate the result but we also indicate very briefly how it is proven (as we will use the technique in another proposition below). 

\prop \label{prop:2.2} 
Suppose that for any element $a$ in $A$ there exist elements $e,f\in A$ so that $ae=a$ and $fa=a$, then $A$ has local units.
\eprop

\bew
i) First it is shown that for any finite number of elements $a_1, a_2, \dots, a_n$ in $A$, there is an element $e\in A$ so that 
$a_ie=a_i$ for all $i$. This is done by induction. Suppose that the result is true for any set of $n-1$ elements. Then first choose $e'$ in $A$ satisfying $a_1e'=a_1$ and then take $e''$ in $A$ so that $(a_i-a_ie')e''=a_i-a_ie'$ for all $i=2,3,\dots,n$. One verifies easily that $e$, defined as $e'+e''-e'e''$, will be the element we are looking for. Similarly, we can find an element $f$ so that $fa_i=a_i$ for all $i$.
\vskip 3pt
ii) Next, take again a finite number of elements $a_1, a_2, \dots, a_n$. By the result above, we find elements $e'$ and $e''$ in $A$ so that 
\begin{equation*}
a_ie'=a_i \qquad\qquad\text{and}\qquad\qquad e''a_i=a_i
\end{equation*}
for all $i=1,2,\dots,n$. Now we define again $e=e'+e''-e'e''$ and it is again straigthforward to verify that this element will satisfy
\begin{equation*}
a_ie=a_i \qquad\qquad\text{and}\qquad\qquad ea_i=a_i
\end{equation*}
for all $i=1,2,\dots,n$ as required.
\ebew

Using the above technique, we can prove the following result for a $^*$-ring. 

\prop
If $A$ is a $^*$-ring with local units as in Definition \ref{defin:2.1}, then we can also find self-adjoint (and even positive) local units.
\eprop

\bew
Take as before a finite number of elements $a_1, a_2, \dots, a_n$ in $A$. As we assume the existence of local units, it is possible to find an element $f\in A$ so that
\begin{equation*}
a_if=a_i \qquad\qquad\text{and}\qquad\qquad a_i^*f=a_i^*
\end{equation*}
for all $i=1,2,\dots,n$. Then also $f^*a_i=a_i$ for all $i$. Now let $e=f +f^*{\joost -}ff^*$. Just as in ii) of the proof of the previous proposition, we now will have
\begin{equation*}
 a_ie=a_i \qquad\qquad\text{and}\qquad\qquad ea_i=a_i
\end{equation*}
for all $i=1,2,\dots,n$. And obviously $e=e^*$. Finally, by replacing $e$ by $e^2$ if necessary, we find that it is possible to get positive local units in $A$.
\ebew

Recall that a ring $A$ is {\em idempotent} if $A^2=A$. A ring is called {\em firm} (also called {\em tensor idempotent}) if the multiplication induces an isomorphism $A\ot_AA\stackrel{\cong}{\longrightarrow} A$. The surjectivity of the latter map means exactly that $A$ is idempotent, hence any firm ring is idempotent, but the converse is not necessarily true. It was however proven in \cite[Theorem 1.1]{BV} that if $A$ is idempotent, then $A\ot_A A$ is firm. Explicitly, this means that the map 
$$ (A\ot_A A)\ot_{A\ot_A A} (A\ot_A A) \to A\ot_A A,\ (a\ot_A a')\ot (b\ot_A b')\mapsto aa'b \ot_A b' = aa'\ot_A bb' = a\ot_A a'bb'$$
is bijetive. Let us now show that examples of firm rings arise from rings with local units.

\lem
Let $A$ be a ring with (left) local units. Then $A$ is (left) non-degenerate, idempotent and firm. 
\elem

\bew
Take $a\in A$. If $ba=0$, for all $b\in A$, then in particular $ea=a=0$, where $e\in A$ is a left local unit for $a$, hence $A$ is left non-degenerate. 

Define $\ol\mu:A\to A\ot_A A$ by $\ol\mu(a)=e\ot_A a$, where $e\in A$ is any left local unit for $A$. Then $\ol \mu$ is indeed a well-defined map, since of any other $e'\in A$ such that $e'a=a$, we have that 
$$e\ot_A a=e''e\ot_A a=e''\ot_A ea=e''\ot_A a=e''\ot_A e'a=e''e'\ot_A a=e'\ot_A a,$$
where $e''$ is a common left local unit for $e$ and $e'$. Obviously, $\ol\mu$ is a right inverse for the map $A\ot_A A\to A$ induced by multiplication. Let us check that it is also a left inverse. Take any $\sum_i a_i\ot_A a'_i\in A\ot_A A$. Let $e$ be a left local unit for all $a_i$. Then $e$ is also a left local unit for $\sum_i a_ia'_i$, and hence 
$$\ol\mu\circ \mu(\sum a_i\ot a'_i)=\sum_i e\ot_A a_ia'_i=\sum_i ea_i\ot_A a'_i=\sum_i a_i\ot a'_i.$$
This proves that $A$ is firm, hence idempotent.
\ebew

Before we continue, let us give 
some simple but typical examples of rings with local units, as well as examples of rings that are non-degenerate and idempotent but do not have local units.
More examples will be given and discussed in Section \ref{s:examples}.

\voorb\label{voorb:2.3}
i) Let $X$ be any set and let $A=K(X)$ be the algebra of complex valued functions with finite support in $X$ (as in Example \ref{voorb:1.2}i)). This algebra clearly has local units. The local units can be chosen to be idempotents. If $X$ is infinite, there is no unit.
\ssnl
ii) Take for $X$ the set $\mathbb N$ of natural numbers in example i) above, as we also did in Example \ref{voorb:1.2}ii). Then elements in $K(X)$ are sequences $(a_n)_{n=1}^\infty$ with $a_n\in \mathbb C$ but so that only finitely many of these elements are non-zero. We can also consider the set $c_0$ of sequences that tend to $0$ when $n\to\infty$. This is also a non-degenerate algebra, but now there are no local units.
Indeed, if e.g.\ $ae=a$ with $a_n\neq 0$ for all $n$, then we must have $e_n=1$ for all $n$ and $1$ does not belong to $c_0$.
\ssnl
iii) Let $A$ be the algebra of continuous maps $\mathbb R\to \mathbb R$ with a compact support. Let $f$ be any such function and let $[a,b]$ be any interval containing the support of $f$. Then define $e:\mathbb R \to \mathbb R$ as the function which is defined as $0$ outside $[a-1,b+1]$, $1$ on $[a,b]$, $x-a+1$ on $[a-1,a]$ and $b-x+1$ on $[b,b+1]$. Then $e$ is a two-sided local unit of $f$. However $e$ nor any other local unit is idempotent (in fact, the only idempotent element of $A$ is $0$). 
\evoorb

The algebra $c_0$ is still idempotent. Here is an example where this is not even true.

\prop\label{prop:2.4}
Let $A$ be the algebra of sequences $(a_n)_{n\in\mathbb N}$ of complex numbers with the property that $(na_n)$ is bounded. This is a subalgebra of $c_0$. It is non-degenerate. There is no unit, there are no local units and the algebra is not idempotent.
\eprop

\bew
i) It is clear that $A$ is a vector space. Now assume that $a,b\in A$. Because $|a_n|\leq |na_n|$ for $n\geq 1$ the sequence $a$ will itself also be bounded. If $|a_n|\leq M$ for all $n$, then $|na_nb_n|\leq M|nb_n|$ for all $n$ and it follows that $ab\in A$. Hence $A$ is an algebra.
\vskip 3pt
ii) Now assume that $a\in A$ and that $ab=0$ for all $b\in A$. Take for $b$ the sequence with $0$ in all places, except $b_k=1$ for a given $k$. Then $a_k=0$. This is true for all $k$ so that $a=0$. Because the algebra is abelian, we only have to consider this case. Hence it is non-degenerate.
\vskip 3pt
iii)
Now we show that it is not idempotent. Indeed, if $a,b\in A$, then $(n^2a_nb_n)$ will be a bounded sequence. It follows that $(n^2c_n)$ is a bounded sequence for all $c\in A^2$. Now consider the sequence $c$  given by $c_n=\frac1n$ for all $n$. Clearly $n^2c_n=n$ for all $n$ and so $c\notin A^2$.
\ebew

We have mentioned already that the multiplier algebra of $c_0$ consist of all bounded sequences (see Example \ref{voorb:2.11c}). Now any bounded sequence is clearly also a multiplier of the subalgebra of $c_0$ described in Proposition~\ref{prop:2.4}. This provides an example of a proper subalgebra that has the same multiplier algebra as the original one.

\nl
\bf The strict topology on $M(A)$ \rm 
\nl
Let $T$ be a ring and $A$ an non-degenerate subring of $T$.
We now define a topology on $T$ in terms of convergence of nets.

\defin
A net $(x_\alpha)$ in $T$ converges to an element $x\in T$ in the {\it left strict topology} if for all $a\in A$ we have that $xa=x_\alpha a$ for $\alpha$ large enough.
\edefin
We have that $x_\alpha\to x$ if and only if for any finite set $\{a_1,a_2,\dots,a_n\}$ of $A$, we have that $xa_i=x_\alpha a_i$ for all $i$ when $\alpha$ is large enough. Similarly, we introduce a right strict topology and a two-sided strict topology (which we call simply the {\em strict topology}). In the latter, we have that A net $(x_\alpha)$ in $T$ converges to an element $x\in T$ if for all $a\in A$ we have that $xa=x_\alpha a$ and $ax=ax_\alpha$ for $\alpha$ large enough.
\snl

The following follows from standard arguments, and has an obvious one-sided version.
\lem
A basis of open sets in the strict topology is given by sets of the form 
$$\mathcal O(x,a_1,\ldots,a_n)=\{y\in T~|~ a_iy=a_ix \text{ and } ya_i=xa_i, i=1,\ldots, n\}$$ indexed by any $x\in T$ and $a_1,\ldots,a_n\in A$. This topology is also known as the {\em $A$-adic} topology on $T$, see e.g.\ \cite[Section 41.4]{BW}.
\elem.

We have the following result (see \cite[Proposition 3.9]{CVW2}, although formulated differently). 

\prop\label{prop:new2.7}
Let $A$ be a non-degenerate ring that is a subring of a unital ring $T$. Then the following
assertions are equivalent:
\begin{enumerate}
\item[(i)] $A$ is dense in $T$ with respect to the (left) strict (or $A$-adic) topology,
\item[(ii)] $A$ is a (left) ideal of $T$ and has (left) local units.
\end{enumerate}
\eprop 

\bew
$\underline{(i)\Rightarrow(ii)}$ If $A$ is dense in $T$, then any open set $\mathcal O(x,a_1,\ldots,a_n)$ has a non-empty intersection with $A$. Hence, for any $x\in T$ there exists $a\in A$ such that $xa_i=aa_i\in A$. Therefore $A$ is an ideal. Taking in particular $x=1$, we find that there exists $e\in A$ such that $a_i=1a_i=ea_i$, so $e$ is a left local unit for the elements $a_i$. The two-sided version is proven in the same way.
\ssnl
$\underline{(ii)\Rightarrow(i)}$  
Since $A$ has left local units, for all $a_1,\cdots,a_n\in A$, we can 
find $e\in A$ such that $ea_i=a_i$. Then every $x\in T$, we 
then have $xea_i=xa_i$. Since $A$ is a left ideal in $T$, we find $xe\in A\cap \mathcal O(x,a_1,\ldots,a_n)$, so $A$ is dense in $T$.
\ebew

\prop
If $A$ is a ring with local units, then for a left $A$-module the following statements are equivalent:
\begin{enumerate}
\item[(i)] $M$ is firm (i.e.\ $A\ot_A M\cong M$ via multiplication);
\item[(ii)] $M$ is idempotent (i.e.\ $AM=M$);
\item[(iii)] $A$ acts with local units on $M$ (i.e.\ for any $m_1,\ldots,m_n\in M$ there exists $e\in A$ such that $em_i=m_i$ for $i=1,\ldots,n$);
\item[(iv)] $A$ is dense in $M(A)$ with respect to the $M$-adic topology, which is the topology generated by the basis of open sets ${\mathcal O}(x,m_1,\ldots,m_n)=\{y\in M(A)~|~xm_i=ym_i, i=1,\ldots,n\}$ for all $x\in M(A)$ and $m_1,\ldots,m_n\in M$.
\end{enumerate}
\eprop

\bew
If $M$ is firm, then the map $A\ot_A M\to M$ induced by multiplication is an isomorphism, hence in particular surjective, which means exactly that $AM=M$.

If $AM= M$, then any $m$ can be written as $\sum_i a_im_i$. Consider now a local unit $e$ for the elements $a_i$, then we also find that $e$ is a local unit for $m$: $em=\sum_i ea_im_i=\sum_i a_im_i=m$.

If a left $A$-module is unital, then it is also firm since the morphism $M\to A\ot_AM, m\mapsto e\ot_A m$ with $e$ a local unit is an inverse for the morphism $A\ot_A M,\ a\ot_A m\mapsto am$.

The equivalence with the last statement is proven in the same way as Proposition \ref{prop:new2.7}.
\ebew

Proposition~\ref{prop:new2.7} can in particular be applied to the case $T=M(A)$.
Sometimes, we call $A$ a dense ideal because it is an essential ideal in $M(A)$. We see now that this is appropriate only if $A$ has local units. Consequently, a natural question that arises is what can we say in general, when $A$ is not assumed to have local units? In this case, we still can consider the closure $\overline A$ of $A$ in $M(A)$ as in the following result.

\prop
Denote by $\overline A$ the space of elements in $M(A)$ that can be approximated by a net in $A$ in the strict topology on $M(A)$.  Then $A$ is a two-sided ideal of $\overline A$ and $\overline A$ in turn is a two-sided ideal of $M(A)$.
\eprop

\bew
The first statement is obvious as $\overline A$ is a subspace of $M(A)$ by definition. 
\vskip 3pt
To prove the second statement, take any $y$ in $M(A)$ and $x\in\overline A$. Take a net $(x_\alpha)$ in $A$ that converges to $x$. Then $(x_\alpha y)$ will be a net in $A$ that converges to $xy$ and $(y x_\alpha)$ will be a net that converges to $yx$. The argument is similar as in the first part of the proof of Proposition \ref{prop:new2.7}..
\ebew


Furthermore, we can show the following characterization of the multiplier algebra.

\prop
$M(A)$ is the completion of $A$ with respect to the strict topology on $A$.
\eprop

\bew
Assume that $(x_\alpha)$ is a Cauchy net in $A$ for the strict topology. This implies that for all $a$ in $A$ we have an index $\alpha_0$ such  that  $a x_\alpha-a  x_\beta=0$ and $x_\alpha a-x_\beta a=0$ for $\alpha,\beta\geq \alpha_0$. Define the maps $a\mapsto ax$ and $a\mapsto xa$ from $A$ to itself by $ax=ax_{\alpha_0}$ and $xa=x_{\alpha_0}a$. Then $ax=ax_\alpha$ and $xa=x_\alpha$ when $\alpha$ is big enough.
\ssnl
For all $a,b\in A$ we have, when $\alpha$ is big enough, 
\begin{equation*}
(ax)b=(ax_\alpha)b=a(x_\alpha b)=a(xb).
\end{equation*}
So we have $x\in M(A)$. 
\ssnl
By definition $xa=\lim_\alpha x_\alpha a$ and $ax=\lim_\alpha ax_\alpha$. 
\ssnl
This proves that $M(A)$ is the completion of $A$ for the strict topology on $A$.
\ebew

It follows that $M(A)$ is complete for the strict topology.

%
%
%

%
%

\section{\hspace{-17pt}. Examples and special cases} \label{s:examples} 

In this section, we will discuss the existence of local units for multiplier Hopf algebras, algebraic quantum groups, algebraic quantum hypergroups and weak multiplier Hopf algebras. 
But we start with a class of examples of another nature, namely the rational dual of a co-Frobenius coalgebra.
\nl
\bf The rational dual of co-Frobenius coalgebras \rm
\nl
Let $C$ be a $k$-coalgebra (over a field $k$), and $C^*=\Hom_k(C,k)$ its dual algebra, whose multiplication is given by the formula
$$(f*g)(c)=\sum_{(c)}f(c_{(1)})g(c_{(2)}),$$
for all $f,g\in C^*$ and $c\in C$. We use the Sweedler notation $\sum_{(c)} c_{(1)}\ot c_{(2)}$ for $\Delta(c)$. 
Similarly, any right $C$-comodule is a left $C^*$-module by means of the action
$$f.m=\sum_{(m)}f(m_{[1]})m_{[0]},$$
for all $f\in C^*$ and $m\in M$.
Recall that the {\em (left) rational} part of a left $C^*$-module $M$ is the following subset
$${\sf Rat}({_{C^*}M})=\{m\in M ~|~ \exists m_i\in M, c_i\in C, i=1,\ldots, n : f*m=\sum_i f(c_i)m_i, \forall f\in C^*\}.$$
Clearly, ${\sf Rat}({_{C^*}M})$ is a left $C^*$-submodule of $M$. In fact, one can show that the left rational part of $C^*$ is the biggest $C^*$-submodule that is induced by right a $C$-comodule. In case ${\sf Rat}(M)=M$, we say that $M$ is a rational $C^*$-module.
In particular, we have
$${\sf Rat}({_{C^*}C^*})=\{f\in C^* ~|~ \exists f_i\in C^*, c_i\in C, i=1,\ldots, n : g*f=\sum_i g(c_i)f_i, \forall g\in C^*\}$$
which is a (rational) left ideal in $C^*$.

Recall (see e.g. \cite[Definition 3.2.4]{DNR}) that $C$ is called {\em right semiperfect} if and only if ${\sf Rat}({_{C^*}C^*})$ is dense in $C^*$ with respect to the finite topology. This finite topology is the topology generated by the basis of open sets 
$${\mathcal O}(f,c_1,\ldots,c_n)=\{g\in C^*~|~g(c_i)=f(c_i), i=1,\ldots,n\}$$
indexed by all $f\in C^*$ and $c_1,\ldots, c_n\in C$.  In the same way, one defines right and two-sided semi-perfect rings. The next result also has a right and two-sided version.

\prop \label{prop:semiperfect}
Let $R$ be a rational left ideal in $C^*$. If $R$ is dense in the finite topology of $C^*$, then $R$ is also dense in the left $R$-adic (i.e. the strict) topology on $C^*$, hence $C^*$ has left local units.
In particular, if $C$ is right semi-perfect, then ${\sf Rat}({_{C^*}C^*})$ has left local units.
\eprop

\bew
Consider $f\in C^*$ and $r^1,\ldots,r^n\in R$. By rationality of $R$, we have $fr^i = \sum_i f(c^i_j)r^i_{j}$ for some $c^i_j\in C$ and $r^i_j\in R$. If $R$ is dense in the finite topology on $C^*$, then there exist an $s\in R$ such that $s(c^i_j)=f(c^i_j)$, hence $fr^i = \sum_j f(c^i_j)r^i_{j}=\sum_j s(c^i_j)r^i_{j}=sr^i$, and $R$ is dense in the $R$-adic topology. The existence of local units then follows from Proposition \ref{prop:new2.7}.
\ebew

One can show (see \cite[Corollary 3.2.17]{DNR} that when $C$ is both left and right semi-perfect then ${\sf Rat}({}_{C^*}C^*)={\sf Rat}(C^*_{C^*})$, which is then by consequence a ring with (two-sided) local units.

The following provides important classes of semi-perfect coalgebras.
A coalgebra $C$ is called right co-Frobenius if and only if there exists an injective right $C^*$-linear map $j:C\to C^*$. More generally, $C$ is called right quasi-co-Frobenius if and only if there exists an injective left $C^*$-linear map $j:C\to (C^*)^I$ for some (possibly infinite) index set $I$.  There are obvious left and two-sided versions of these notions. For other characterizations of (quasi-)co-Frobenius coalgbras, we refer to \cite[Section 3.3]{DNR}. Let us just remark that a finite dimensional coalgebra is (quasi-)co-Frobenius if and only if its dual algebra is (quasi-)Frobenius. The next result combines some known results with the earlier results of this paper.

\prop
Let $C$ be a right quasi-co-Frobenius coalgebra, with Frobenius morphism $j:C\to (C^*)^I$.
i) The morphism
$$\tilde\jmath: C^{(I)}\to C^*, \tilde\jmath(c)(d)=j(d)(c)$$
is left $C^*$-linear and its image is contained in ${\sf Rat}(_{C^*}C^*)$ and is dense in $C^*$ with respect to the finite topology.  Consequently, $C$ is right semiperfect and hence ${\sf Rat}(_{C^*}C^*)$ has left local units. 
(Here we denoted by $C^{(I)}$ a direct sum of $I$ copies of $C$.)
\ssnl
ii) If $C$ is right co-Frobenius, then ${\sf Im} j\subset {\sf Rat}(C^*_{C^*})$ is an essential ideal right ideal in $C^*$ and hence $C^*$ can be interpreted as the multiplier algebra of ${\sf Rat}(C^*_{C^*})$.
\ssnl
iii) If $C$ is a two-sided co-Frobenius algebra, then ${\sf Rat}(_{C^*}C^*)={\sf Rat}(C^*_{C^*})$ is an algebra with local units whose multiplier algebra is isomorphic to $C^*$.
\eprop

\bew
i) 
Let us check that $\tilde\jmath$ is left $C^*$-linear. Consider any $c\in C^{(I)}$ and $d\in C$. Then we find indeed that
\begin{eqnarray*}
\tilde\jmath(f.c)(d) &=&\textstyle\sum_{(c)} \tilde\jmath(f(c_{(2)})c_1)(d) = \textstyle\sum_{(c)}f(c_{(2)})\tilde\jmath(c_1)(d) \\
&=&\textstyle\sum_{(c)} f(c_{(2)})j(d)(c_1) = (j(d)*f)(c) = j(d.f)(c)\\
&=& \tilde\jmath(c)(d.f) = \textstyle\sum_{(d)}\tilde\jmath(c)(f(d_{(1)})d_{(2)}) \\
&=& \textstyle\sum_{(d)} f(d_{(1)})\tilde\jmath(c)(d_{(2)}) = (f*\tilde\jmath(c))(d)
\end{eqnarray*}
The computations above also show that $f*\tilde\jmath(c)=\sum_{(c)}f(c_{(2)})\tilde\jmath(c_1)$, which means that ${\sf Im}\tilde\jmath$ is contained in ${\sf Rat}(_{C^*}C^*)$. 
\ssnl
Now suppose that $c\in C$ is such that $\tilde\jmath(d)(c)=0$ for all $d\in C^{(I)}$. Then we find that $j(c)(d)=0$ for all $d\in C$, hence $j(c)=0$. Since $j$ is injective, this means that $c=0$. This shows that ${\sf Im}\tilde\jmath^\bot=0$, hence ${\sf Im}\tilde\jmath$ is dense with respect to the finite topology on $C^*$ by \cite[Theorem 1.2.6]{DNR}. Since ${\sf Im}\tilde\jmath$ is contained in ${\sf Rat}(_{C^*}C^*)$, the latter is also dense in $C^*$ and therefore $C$ is right semiperfect.
\ssnl
ii) From similar computation as in part $1$, it follows that ${\sf Im} j\subset {\sf Rat}(C^*_{C^*})$.
Suppose that $f\in C^*$ is such that $g*f=0$ for all $g\in {\sf Rat}(_{C^*}C^*)$. Then in particular, we have that
for all $c\in C$: $0=j(c)*f=j(c.f)$. Since $j$ is injective, this implies that $0=c.f=f(c_{(2)})c_{(1)}$. When we apply $\epsilon$ to the last equation, we find that $f(c)=0$, for all $c\in C$, hence $f=0$. This shows that ${\sf Im j}$ is an essential ideal in $C^*$. The last statement then follows from Proposition \ref{prop:Misbiggest}.
\ssnl
iii) If $C$ is left and right co-Frobenius, then $C$ is left and right semi-perfect by part i). Therefore the left and right rational dual coincide. The remaining assertions then follow directly from part ii) and Proposition \ref{prop:semiperfect}.
\ebew

The results above can be generalized to corings, see \cite{CI}, \cite{IV}.
The results of this section lead to the following classes of algebras with local units.
\ssnl
If $X$ is any set, and $kX$ is the coalgebra defined by considering all elements of $X$ as grouplike, that is, $\Delta(x)=x\ot x$ for all $x\in X$. Then the algebra $kX^*$ is the algebra of all $k$-valued functions on $X$, with pointwise operations. Denote as before by $\delta_x\in kX^*$ the function that sends $x$ to $1$ and all other elements of $X$ to $0$. Then the $k$-linear map $j:kX\to kX^*, x\mapsto \delta_x$ is $kX^*$-linear since
$$j(f.x) = j(f(x) x)=f(x)\delta_x = f*\delta_x=f*j(x).$$ 
Hence $kX$ is co-Frobenius. The rational dual of $kX$ is exactly the algebra $k(X)$ of $k$-valued functions with finite support, that we encountered before in Examples \ref{voorb:1.2} and \ref{voorb:1.2bis}.
\ssnl
If $H$ is not just a coalgebra, but a Hopf algebra, then $H$ is left or equivalently right quasi-co-Frobenius if and only if it is left or equivalently right co-Frobenius if and only if there exists a non-zero (left or right) integral in $H^*$. In this case the rational dual of $H$ is a multiplier Hopf algebra. 
\ssnl
Our next aim is to show that in fact any multiplier Hopf algebra has local units.

\nl
\bf Multiplier Hopf algebras \rm
\nl
It was shown first in \cite{Dr-VD-Z} that any regular multiplier Hopf algebra has one-sided local units. More precisely, in Proposition 2.2 of \cite{Dr-VD-Z} it is proven that for any finite number of elements $a_1,a_2,\dots,a_n$ in a regular multiplier Hopf algebra $A$, there exists elements $e,f$ so that $ ea_i=a_i$ and $a_i f=a_i$ for all $i$. In the proof it is used that the antipode $S$ maps $A$ into $A$. In fact, the result with the same proof is still valid for any multiplier Hopf algebra (regular or not)  with this property. Further, in Proposition 2.6 of the same paper, it is shown that we can assume $e=f$ in the case of an algbraic quantum group (i.e.\ a regular multiplier Hopf algebra with integrals).
\snl
Now, due to the result recalled in Proposition \ref{prop:2.2} of this paper, we know that it already follows that $A$ has local units whenever we have a multiplier Hopf algebra with an antipode that maps $A$ to itself. There is no need to have integrals.
\snl
Next, in \cite{VD-Z}, it is proven that for any multiplier Hopf algebra $A$, given $a$ in $A$, there exists elements $e$ and $f$ in $A$ so that $ea=a$ and $af=a$ (see Proposition 1.2 in \cite{VD-Z}). Again, because of the result in Proposition \ref{prop:2.2} it follows that $A$ has local units.
\snl
Therefore, we have the following result. We recall also the proof of the result in \cite{VD-Z}, not only for the convenience of the reader, but also because we will need similar techniques in later results.

\prop\label{prop:3.1}
Let $(A,\Delta)$ be a multiplier Hopf algebra. Then $A$ has local units.
\eprop
\bew
i) First, we claim that for any $a,b\in A$ we have that $\Delta(b)(1\ot a)\in A\ot Aa$. Indeed, suppose that $\omega$ is a linear functional on $A$ so that $\omega(qa)=0$ for all $q$ in $A$. Then 
\begin{equation*}
(\iota\ot\omega)((c\ot 1)\Delta(b)(1\ot a))=0
\end{equation*}
for all $c$. As the product in $A$ is non-degenerate, and because $\Delta(b)(1\ot a)\in A\ot A$, it follows that also $(\iota\ot\omega)(\Delta(b)(1\ot a))=0$. This implies that $\Delta(b)(1\ot a)\in A\ot Aa$.
\ssnl
ii) Now write $\Delta(b)(1\ot a)=\sum_i p_i\ot q_ia$. Then $\sum_i S(p_i)q_i a= \varepsilon(b)a$. If now we take $b$ so that $\varepsilon(b)=1$, we find an element $e$ in $A$ so that $ea=a$. 
\ssnl
iii) Similarly we get that $(a\ot 1)\Delta(b)\in aA\ot A$ for all $a,b$ and with $\varepsilon(b)=1$ we again will get an element $f$ in $A$ so that $af=a$. Then, using Proposition \ref{prop:2.2} we find that $A$ has local units.
\ebew

The argument can be refined in such a way that we directly obtain local units, without using the result of Proposition \ref{prop:2.2}, but of course this is more complicated. 
\nl
\bf Algebraic quantum hypergroups \rm
\nl
The existence of local units for algebraic quantum hypergroups is proven in Proposition 1.6 of \cite{De-VD}. The proof is inspired by the one given in \cite{Dr-VD-Z}.  We give here another argument, partly motivated by the one given above, in the proof of Proposition \ref{prop:3.1} for multiplier Hopf algebras. 

\prop\label{prop:3.2}
Let $(A,\Delta)$ be an algebraic quantum hypergroup. Then $A$ has local units.
\eprop

\bew
i) To prove that $\Delta(b)(1\ot a)\in A\ot Aa$ we need that 
\begin{equation*}
\Delta(A)(1\ot A)\subseteq A\ot A
\tussenen
(A\ot 1)\Delta(A)\subseteq A\ot A
\end{equation*}
and that the product in $A$ is non-degenerate. These properties are still valid, see Definition 1.1 in \cite{De-VD}. So we can proceed here as in the proof of Proposition \ref{prop:3.1} above.
\ssnl
ii) For the second step, we need another argument because the antipode does not satisfy the formula used in the previous proof. Instead, we can  use the right integral $\psi$. Indeed, we have $\Delta(b)(1\ot a)\in A\ot Aa$ and if we apply $\psi$ on the first factor, we obtain $\psi(b)a\in Aa$. By choosing $b$ such that $\psi(b)=1$ we get $a\in Aa$. Therefore we have an element $e\in A$ satisfying $ea=a$. 
\ssnl
iii) In a similar way, we have an element $f$ satisfying $af=a$.  Then we apply Proposition \ref{prop:2.2} to complete the proof.
\ebew

The argument in the proof above can also be used to show that algebraic quantum groups admit local units. In fact, as any algebraic quantum group is an algebraic quantum hypergroup, the result would follow from the above proposition.
\nl
\bf Weak multiplier Hopf algebras \rm
\nl
We consider  a weak multiplier Hopf algebra $(A,\Delta)$ as defined and studied in \cite{VD-W0}. We use $\varepsilon_s$ for the source map $a\mapsto \sum_{(a)}S(a_{(1)})a_{(2)}$ and $\varepsilon_t$ for the target map $a\mapsto \sum_{(a)}a_{(1)}S(a_{(2)})$. The range of the source map, the source algebra, is denoted by $B$ while the range of the target map, the target algebra by $C$. See Notation 2.9 in \cite{VD-W1}. 
\ssnl
To show that the underlying algebra of a weak multiplier Hopf algebra, or of a algebraic quantum groupoid, has local units is more involved.
The proof is given in Proposition 2.21 of \cite{VD-W1}. We will not repeat the full proof here, we just give some ideas involved in the different steps. 
\ssnl
The first part, as in the following lemma, is of the same nature as in the proof of the existence of local units for multiplier Hopf algebras in Proposition \ref{prop:3.1}.

\lem\label{lem:3.3}
For all $a\in A$ we have that $Ca\subseteq Aa$ and $aB\subseteq aA$.
\elem

\bew
i) Fix $a\in A$. Assume that $\omega$ is a linear functional such that $\omega(ba)=0$ for all $b\in A$. Then
\begin{equation*}
(\iota\ot\omega)((1\ot b)((\iota\ot S)((c\ot 1)\Delta(p)))(1\ot a))=0.
\end{equation*}
for all $b,c,p\in A$. We use that $(c\ot 1)\Delta(p)\in A\ot A$ and that $S$ maps $A$ into $M(A)$. We also know that $((\iota\ot S)\Delta(p))(1\ot a)$ belongs to $A\ot A$ (see Remark 1.2 in \cite{VD-W1}. Then we can cancel $c$ and obtain that 
\begin{equation*}
(\iota\ot\omega)((1\ot b)((\iota\ot S)(\Delta(p))(1\ot a))=0.
\end{equation*}
We write $((\iota\ot S)(\Delta(p)))(1\ot a)=\sum_i p_i\ot q_i$ and assume that the elements $p_i$ are linearly independent. So we find $\omega(bq_i)=0$ for all $i$ and all $b$. Replace $b$ by $p_i$ and take the sum over $i$. Then
\begin{equation*}
\omega(\varepsilon_t(p)a)=\sum_i \omega(p_iq_i)=0.
\end{equation*}
We see that $\omega(xa)=0$ for all $x$ in the target algebra $C$.
\ssnl
ii) Again take $a\in A$ and now assume that $\omega$ is a linear functional on $A$ such that $\omega(ab)=0$ for all $b\in A$. Then \begin{equation*}
(\omega\ot\iota)((a\ot 1)(S\ot\iota)(\Delta(p)(1\ot c))(b\ot 1))=0.
\end{equation*}
Here we use that $\Delta(p)(1\ot c)\in A\ot A$ and that $S$ maps $A$ to $M(A)$. On the other hand, because $(a\ot 1)(S\ot\iota)(\Delta(p))\in A\ot A$, we can cancel $c$ and we find
\begin{equation*}
(\omega\ot\iota)((a\ot 1)(S\ot\iota)(\Delta(p))(b\ot 1))=0.
\end{equation*}
Now we write $(a\ot 1)(S\ot\iota)(\Delta(p))=\sum_i p_i\ot q_i$ with the $q_i$ linearly independent. Then we find $\omega (p_ib)=0$ for all $i$ and all $b$. Finally, replace $b$ by $q_i$ and take the sum to get
\begin{equation*}
\omega(a\varepsilon_s(p))=\sum_i\omega(p_iq_i)=0.
\end{equation*}
Hence $\omega(ay)=0$ for all $y\in B$.
\ebew

For the next step, we need two properties of the source and target algebras. For the target algebra $C$ we have $CA=A$, see Proposition 2.9 in \cite{VD-W1}. We also know that $C$ has local units, see 
the proof of Proposition 2.14 in \cite{VD-W1}.
\ssnl
It would take us too far to include the argument for the property $CA=A$, but we can sketch the proof of the existence of local units for $C$. We use the canonical idempotent $E$, which belongs to the multiplier algebra $M(B\ot C)$, see Section 2 of \cite{VD-W1}.

\prop\label{prop:3.4}
The source algebra $B$ and the target algebra $C$ have local units.
\eprop

\bew
i) Take an element $u$ in $C$ and assume that $\omega$ is a linear functional on $C$ so that $\omega(xu)=0$ for all $x\in C$. It is shown that $(y\ot 1)E\in B\ot C$ for all $y\in B$. Then
\begin{equation*}
(\iota\ot\omega)((y\ot 1)E(1\ot u))=0
\end{equation*}
for all $y\in B$. Now we use that also $E(1\ot u)\in B\ot C$. Then, using that the product in $B$ is non-degenerate, we get
\begin{equation}
(\iota\ot\omega)(E(1\ot u))=0.\label{eqn:3.1}
\end{equation}
\ssnl
ii) Now we can apply the distinguished functional $\varphi_B$, defined and characterized by the property that $(\varphi_B\ot\iota)E=1$. This gives $\omega(u)=0$. It proves that $u\in Cu$ and we get an element $e\in C$ satisfying $u=eu$.
\ssnl
iii) Similarly we get $u\in uC$ . It follows that $C$ has local units.
\ssnl
iv) The same technique is used to show that $B$ has local units.
\ebew

The element $E$ is a separability idempotent in the multiplier algebra $M(B\ot C)$ (as studied in \cite{VD-si.v1, VD-si.v2}). It is shown in Proposition 1.10 of \cite{VD-si.v2} that the underlying algebras $B$ and $C$ of any separability idempotent have local units. The proof is as above. Another argument uses that $S(E_{(1)})E_{(2)}=1$ where $E_{(1)}\ot E_{(2)}$ is a Sweedler type notation for $E$. Then we apply $m(S\ot\iota)$ (where $m$ is multiplication in $C$) to Equation (\ref{eqn:3.1} and we also arrive at $\omega(u)=0$.
\ssnl
Finally we combine the two results in the existence proof below.

\prop\label{prop:3.5}
Let $(A,\Delta)$ be a weak multiplier Hopf algebra. Then $A$ has local units.
\eprop 

\bew
Take $a\in A$.   Because we know that $a\in CA$, we can write $a=\sum_i a_ix_i$ where $a_i\in A$ and $x_i\in C$. Then we use that $C$ has local units and we find $e\in C$ satisfying 
\begin{equation*}
ae=\sum_i a_ix_ie=\sum_i a_ix_i=a.
\end{equation*}
Next we use the lemma and so $ae\in aA$. Hence $a\in aA$. Then we find an element $f\in A$ satisfying $a=af$. Similarly on the other side and we can conclude that there are local units in $A$.
\ebew

\opm
It is interesting to compare the proofs of Proposition \ref{prop:3.2}, Lemma \ref{lem:3.3} and Proposition \ref{prop:3.4}. Among other things, there are the two possibilites one can use, either with the application of an integral, or with a property of the antipode. 
\eopm

%
%

\section{\hspace{-17pt}. Conclusions} \label{s:conclusions} 

In this note we have collected some material about multiplier algebras and local units.
\ssnl
We discussed the notions of left and right multipliers $L(A)$ and $R(A)$ of a non-degenerate algebra $A$. The multiplier algebra $M(A)$ can be viewed as the intersection of these two, when interpreted correctly. We gave a couple of simple but illustrative examples. We consider non-degenerate homomorphisms from an algebra $A$ to an algebra $B$ as a particular class of homomorphism $\gamma: A\mapsto M(B)$. These non-degenerate homomorphisms can be uniquely extended to a unital homomorphism from $M(A)$ to $M(B)$. We mentioned what could be done for one-sided multipliers. 
Some more research with specific examples would be welcome.
\ssnl
To illustrate some of these properties, we have included an example of a subalgebra $A$ of the algebra $M_3(\mathbb C)$ of $3\times 3$-matrices over the complex numbers. It is a non-unital, non-degenerate idempotent algebra. We have non-trivial left multipliers and right multipliers, while $M(A)$ contains only $A$ itself and the identity. There are no local units for this algebra as then it would have to be unital itself. 
\ssnl
It would be interesting to look for more examples of such type to illustrate the various properties like the existence of local units, non-degeneracy of the product, etc. We refer to \cite{VD-Ve} where we look at this for infinite matrix algebras.
\ssnl
The existence of local units is shown to be equivalent with the density of $A$ in $M(A)$ when $M(A)$ is endowed with the \emph{strict topology}. In the last section we reviewed several classes of examples of algebras with local units, arising from semi-perfect coalgebras, multiplier Hopf algebras, algebraic quantum hyper groups and weak multiplier Hopf algebras.

%
%

\end{document}